\documentclass{amsart}

\usepackage{psfrag, graphicx, wrapfig}

\newtheorem{theorem}{Theorem}
\newtheorem{lemma}[theorem]{Lemma}
\newtheorem{proposition}[theorem]{Proposition}
\newtheorem{corollary}[theorem]{Corollary}

\theoremstyle{definition}
\newtheorem{definition}[theorem]{Definition}

\theoremstyle{remark}
\newtheorem{remark}[theorem]{Remark}

\numberwithin{equation}{section}

\def\bkE{{\rm I\kern-.22em E}}
\def\E{\bkE} 
\def\bkH{{\rm I\kern-.22em H}}
\def\H{\bkH} 
\def\bkR{{\rm I\kern-.17em R}}
\def\RR{\bkR}
\def\SS{\mathbb{S}}
\def\t{a}
\def\e{b}
\def\tri{\mathcal{T}}
\def\G{\mathcal{G}}
\def\K{\kappa}
\DeclareMathOperator{\im}{im}
\newcommand{\co}{\colon\thinspace} 

\begin{document}

\title{Angle structures and normal surfaces} 

\author{Feng Luo}
\address{Department of Mathematics, Rutgers University, New Brunswick, NJ 08854, USA}
\email{fluo@math.rutgers.edu}

\author{Stephan Tillmann}
\address{D\'epartment de math\'ematiques,
  Universit\'e du Qu\'ebec \`a Montr\'eal,
  Case postale 8888, Succursale Centre-Ville,
  Montr\'eal (Qu\'ebec) Canada H3C 3P8} 
\email{tillmann@math.uqam.ca} 




\begin{abstract}
Let $M$ be the interior of a compact 3--manifold with non--empty boundary, and $\tri$ be an ideal (topological) triangulation of $M.$ This paper describes necessary and sufficient conditions for the existence of angle structures, semi--angle structures and generalised angle structures on $(M; \tri)$ respectively in terms of a generalised Euler characteristic function on the solution space of normal surface theory of $(M; \tri).$ This extends previous work of Kang and Rubinstein, and is itself generalised to a more general setting for 3--dimensional pseudo-manifolds.
\end{abstract}

\maketitle


\section{Introduction}
\label{sec:intro}

Let $M$ be the interior of a compact 3--manifold with non--empty boundary, and $\tri$ be an ideal (topological) triangulation of $M.$ Such a triangulation always exists and is not unique. This paper addresses the question whether angle structures on the pair $(M; \tri)$ exist. The main results give necessary and sufficient conditions on the existence of an angle structure, a semi--angle structure and a generalised angle structure respectively in terms of a generalised Euler characteristic function on the solution space of normal surface theory. Our main results are based on and generalise work of Kang and Rubinstein \cite{KR2}.

Angle structures have first been studied by Casson (see \cite{l}) and Rivin \cite{Ri}, and provide a linear analogue of ideal hyperbolic triangulations. An \emph{angle structure} of an ideal triangulation is an assignment of a real number, termed angle, in the range $(0, \pi)$ to each edge of each ideal tetrahedron such that the sum of angles at each vertex of each tetrahedron is $\pi,$ and such that around an edge of $\tri$ the sum of angles is $2\pi.$ Casson and Rivin  observed that the existence of an angle structure implies that all boundary components of $M$ are tori or Klein bottles, and that $M$ is irreducible and atoroidal. For example, if $M$ admits a complete hyperbolic structure of finite volume and $\tri$ is an ideal hyperbolic triangulation of $M,$ then its dihedral angles define an angle structure. In general, given a cusped hyperbolic 3--manifold $M$ of finite volume, subdividing its Epstein--Penner fundamental polyhedron yields a partially flat ideal hyperbolic triangulation of $M,$ where some dihedral angles are $0$ or $\pi.$ This leads to the following two natural extensions of the set of angle structures on $(M; \tri),$ denoted by $A(M; \tri).$

The first considers assignments of angles in the closed interval $[0, \pi],$ called \emph{semi--angle structures} in \cite{KR2}, and the set of all such possible assignments is denoted by $SA(M; \tri).$ If in a semi-angle structure all angles are $0$ or $\pi,$ then it is a \emph{taut} structure as studied by Lackenby \cite{l}. The second extension allows angles to be any real number, giving the set of \emph{generalised angle structures} $GA(M; \tri).$

\begin{theorem}\label{thm: generalised angle structures}
  Let $M$ be the interior of a compact 3--manifold with $v$ boundary
  components, and $\tri$ be an ideal triangulation of $M$ with $t$ ideal
  tetrahedra. Then the following are equivalent:
\begin{enumerate}
\item $(M; \tri)$ admits a generalised angle structure,
\item each boundary component is a torus or a Klein bottle.
\end{enumerate}
If $GA(M; \tri)$ is non--empty, then it is an affine space of
dimension $t+v.$
\end{theorem}

The above theorem was also known to Casson \cite{Ca}. It is, together with Theorems \ref{thm: semi--angle structures} and \ref{thm: angle structures} below, generalised in Section~\ref{sec:pseudo--manifolds} to 3--dimensional pseudo--manifolds with angle sums at vertices and around edges allowing arbitrary real numbers. This provides a natural setting for the study of angle structures on closed 3--manifolds and on 3--orbifolds.

It turns out that the existence of an angle structure is closely related to the theory of normal surfaces.  Summaries of (two different kinds of) normal surface theory can be found in \cite{JR, KR1, ti}. The theory sufficient for this paper is outlined in Section \ref{sec:preliminaries} and uses triangle and quadrilateral coordinates which have to satisfy certain \emph{compatibility equations}. The solution space of the compatibility equations over the real numbers is denoted by $C(M; \tri).$ Work of Kang and Rubinstein shows that $C(M; \tri)$ is a vector space of dimension $t+n,$ where $t$ is the number of ideal 3--simplices and $n$ is the number of ideal 1--simplices in $\tri.$ A canonical basis is given by the so--called \emph{tetrahedral solutions} $W_{\sigma^3_1}, ..., W_{\sigma^3_t}$ and the \emph{edge solutions} $W_{e_1}, ..., W_{e_n}.$ There is a well--defined linear function $\chi^*$ on $C(M; \tri)$ which agrees with the Euler characteristic on embedded and immersed normal surfaces, and gives an upper bound for the Euler characteristic of a branched normal surface. The intersection of $C(M; \tri)$ with the positive unit simplex $\{(x_1, ...,x_{7t}) \in \RR^{7t} | x_i \geq 0, \sum_{i=1}^{7t} x_i =1$\} is a convex rational polytope of dimension $t+n-1,$ called the projective solution space $PC(M; \tri).$  Thus, each point in $PC(M; \tri)$ is a convex linear combination of finitely many \emph{vertex solutions} $V_1,...,V_m.$

\begin{theorem}\label{thm: semi--angle structures}
  Let $M$ be the interior of a compact 3--manifold with non--empty boundary
  and ideal triangulation $\tri.$ Assume that each boundary component
  is a torus or a Klein bottle. Then the following are equivalent:
\begin{enumerate}
\item $(M; \tri)$ admits a semi--angle structure,
\item for all $s \in C(M; \tri)$ with all quadrilateral coordinates
  non--negative, one has $\chi^*(s)\le 0,$
\item for all $s \in PC(M; \tri),$ $\chi^*(s)\le 0,$
\item $\chi^*(V_i)\le 0$ for all vertex solutions $V_i,$ $i=1,...,m.$
\end{enumerate}
\end{theorem}

\begin{theorem}\label{thm: angle structures}
  Let $M$ be the interior of a compact 3--manifold with non--empty boundary
  and ideal triangulation $\tri.$ Assume that each boundary component
  is a torus or a Klein bottle. Then the following are
  equivalent:
  \begin{enumerate}
  \item $(M; \tri)$ admits an angle structure,
  \item for all $s \in C(M; \tri)$ with all quadrilateral coordinates
    non--negative and at least one quadrilateral coordinate positive, $\chi^*(s) < 0,$
  \item for all $s \in PC(M; \tri)$ with at least one quadrilateral coordinate positive, $\chi^*(s) < 0,$
  \item $\chi^*(V_i) < 0$ for all vertex solutions $V_i$ with at least one quadrilateral coordinate positive.
\end{enumerate}
\end{theorem}

\begin{corollary}\label{cor: dimension of angle structures}
  If $(M; \tri)$ admits an angle structure, then the set of all angle
  structures is the interior of a convex rational polytope of
  dimension equal to $t+v,$ where $v$ is the number of ideal vertices, and $t$ is the number of ideal tetrahedra in $\tri.$
\end{corollary}

\begin{remark}
  The main result of \cite{KR2} assumes that $SA(M; \tri)$ is
  non--empty and gives necessary and sufficient conditions for the
  existence of an angle structure; it can be deduced from Theorems
  \ref{thm: semi--angle structures} and \ref{thm: angle structures}
  above. Such ``deformation'' results are interesting because Lackenby
  \cite{l} proves the existence of taut ideal triangulations for an
  interesting class of 3--manifolds. Many of the key ideas in
  \cite{KR2} lead the authors to discover the above results.
\end{remark}


\subsection*{Outline of the paper}

Background material on pseudo-manifolds, ideal triangulations and normal surfaces theory is summarised in Section \ref{sec:preliminaries}, and basic properties of angle structures are given. The theorems stated in the introduction are proved in Section \ref{sec:proofs}. The generalisations of Theorems \ref{thm: generalised angle structures}, \ref{thm: semi--angle structures} and \ref{thm: angle structures} are given in Section \ref{sec:pseudo--manifolds}.


\subsection*{Acknowledgements}

The first author would like to thank Tao Li for discussions and Yair Minksy for hosting his visit to Yale University where part of the work was done. The second author was supported by a Postdoctoral Fellowship from the Centre de recherches math\'ematiques (CRM) and the Institut des sciences math\'ematiques (ISM) in Montr\'eal.


\section{Preliminaries}
\label{sec:preliminaries}

Following \cite{JR}, basic results on triangulations of pseudo-manifolds, ideal triangulations and normal surface theory are recalled in this section. In the final subsection, we give a new (and equivalent) definition of angle structures using normal quadrilateral types.


\subsection{Triangulations of pseudo-manifolds}

Let $ X$ be a union of disjoint compact 3--simplices $\sigma^3_1, ..., \sigma^3_t$ and  $\Phi$ be a collection of affine isomorphisms $\{\phi_1, ...,
\phi_r\}$ such that
\begin{itemize}
\item[(1)] for each $\phi_i,$ there are two distinct codimension-one faces $\tau_i$ and
$\delta_i$ of $\sigma^3_1, ..., \sigma^3_t$ for which $\phi_i:\tau_i \to \delta_i$ is an affine isomorphism, and
\item[(2)] $\{ \tau_i, \delta_i\} \cap \{\tau_j ,\delta_j\} =\emptyset$ for $i \neq j.$ 
\end{itemize}
The quotient space obtained from $X$ by identifying $x \in \tau_i $ with $\phi_i(x) \in \delta_i$ for each $i,$ denoted by $K =X/\Phi,$ is called a \emph{compact
3--dimensional pseudo-manifold}, and it has a cell decomposition $\tri,$ where the cells in $\tri$ are the images of simplices $\sigma^3_1, ..., \sigma^3_t$ and their sub-simplices. The cell decomposition $\tri$ is called a triangulation of $K$ and the cells in $\tri$ are called vertices, edges, triangles and tetrahedra. Let $\tri^{(i)}$ be the set of all $i$-dimensional cells in $\tri.$ If some 2--simplex in $X$ is not identified with another 2--simplex, then $K$ is a \emph{pseudo-manifold with boundary}. Otherwise, it is a \emph{closed pseudo-manifold}. An open $k$--simplex in $X$ embeds into $K$ if $k\in\{0,2,3\}$. The first barycentric subdivision of $(K; \tri)$ is a \emph{$\Delta$--complex} (in the terminology of Hatcher \cite{h}) since now each open simplex embeds into $K$ and its closure has all its vertices distinct. The second barycentric subdivision of $(K; \tri)$ is a simplicial complex. Thus, the potential non-manifold points in $K$ are the vertices and the barycentres of edges in $\tri.$

Suppose a closed pseudo-manifold $K=X /\Phi$ has non-manifold points only at vertices. Depending on the context, one considers a subset $M$ of $K$ which is obtained by deleting a (possibly empty) subset $V'$ consisting of 0--simplices in $K$ or their small open regular neighbourhoods. This gives a pair $(M; \tri'),$ where $\tri'=\{ \sigma \cap M | \sigma \in \tri\}.$ For $\sigma \in \tri,$ the cells $\sigma \cap M$ are called the edges, triangles, and tetrahedra in $\tri'.$ Interesting cases are where $V'=\emptyset$ and $(M; \tri')$ is a triangulated compact 3--manifold, or $V'$ is a small open neighbourhood of $\tri^{(0)}$ and $(M; \tri')$ is a compact 3--manifold with hyperideal triangulation (i.e.\thinspace a partition into truncated tetrahedra), or $V'=\tri^{(0)}$ and $(M; \tri')$ is a non--compact 3--manifold with ideal triangulation. 


\subsection{Normal arcs and discs}

\begin{figure}[t]
\begin{center}
  \includegraphics[width=9cm]{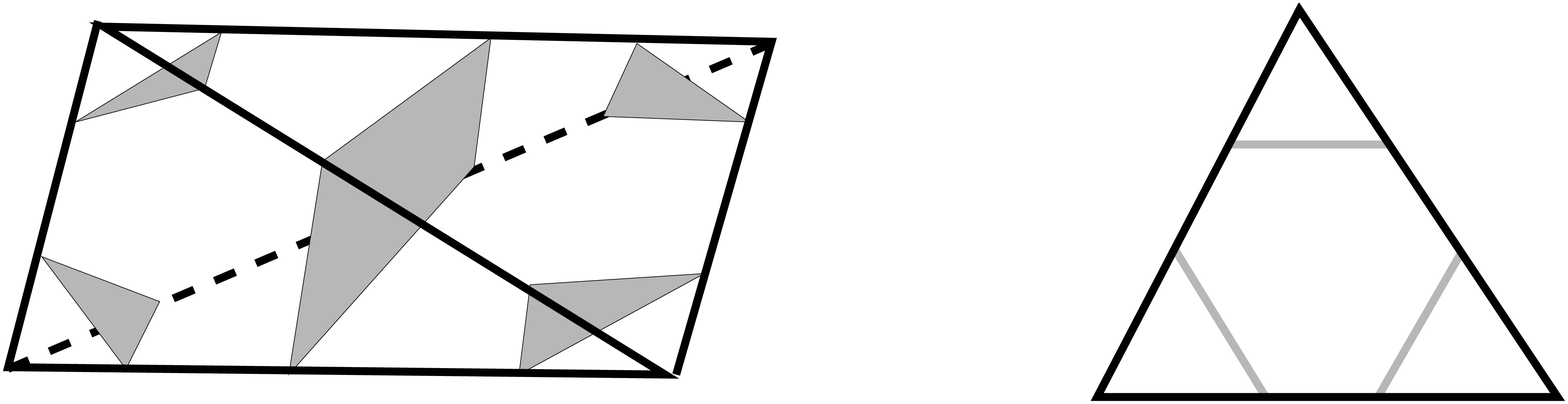}
\end{center}
    \caption{Normal discs in a tetrahedron and normal arcs in a triangle}
    \label{fig:normal discs}
\end{figure}

Given a 2--simplex $\sigma^2,$ a normal arc in $\sigma^2$ is a properly embedded arc in $\sigma^2$ so that its end points are in the interiors of two distinct edges on $\sigma^2.$ For a 3--simplex $\sigma^3,$ a \it normal disc \rm is a properly embedded disc $D$ in $\sigma^3$ so that for each codimension-1 face $\sigma^2$ of $\sigma^3,$ the intersection $D \cap \sigma^2$ is either a single normal arc in $\sigma^2$ or is the empty set. If $\tri$ is a triangulation of a pseudo-manifold $K,$ a \it normal isotopy \rm of $(K, T)$ is an isotopy of $K$ leaving each cell in $\tri$ invariant. There are only seven normal discs in a 3--simplex up to normal isotopy. Four of them are normal triangles and three of them are normal quadrilaterals. Some normal discs are shown in Figure \ref{fig:normal discs}. A normal triangle cuts off a corner from a 3--simplex, and a normal quadrilateral separates a pair of opposite edges. Normal isotopy classes of normal discs are called \emph{normal disc types}. The normal isotopy classes of the quotients of the normal discs in $X$ are called normal disc types in $(K; \tri).$ \emph{Normal arc types} are defined analogously.

\subsection{Normal surfaces}

For a pseudo-manifold $K=X /\Phi$ with a triangulation $\tri,$ a compact surface $F$ in $K$ is called \emph{normal} with respect to $\tri$ if the pre-image of $F$ in $X$ is a collection of normal discs. That is, a normal surface $F$ in $X/\Phi$ is the quotient of a finite collection $Y$ of normal discs in $X$ so that each normal arc $e$ in $Y$ which is not in the boundary of $X/\Phi$ is identified by some $\phi_i$ with another normal arc in $Y.$ A normal surface is \emph{embedded} if the associated collection of normal discs in $X$ is disjoint; it follows from the definition that it is properly embedded.

Fix an ordering of all normal disc types $(q_1, ..., q_{3t}, p_1, ..., p_{4t})$ in $\tri$ where $q_i$ denotes a normal quadrilateral type and $p_j$ a normal triangle type.
The \it normal coordinate \rm $\overline{F} =(x_1,..., x_{3t}, y_1, ..., y_{4t})$ of the normal surface $F$ is
a vector in $\RR^{7t}$ where $x_i$ is the number of normal discs of  type $q_i$ in $F,$ and $y_j$ is the number of normal discs of type $p_j$ in $F.$  An embedded normal surface $F$ is uniquely determined up to normal isotopy by its normal coordinate (see \cite{ti}). We call $(x_1, ..., x_{3t})$ the \emph{quadrilateral coordinate} of $F$ and $\overline{F}.$ 


\subsection{Compatibility equations}

The normal coordinate $\overline{F}$ of a normal surface satisfies a set of homogeneous linear equations, called the \emph{compatibility equations}. They are defined as follows. Suppose $\alpha$ is a normal isotopy class of a normal arc in a triangle $\sigma^2$ in the triangulation $\tri$ which is not in the boundary of $K.$ Then $\sigma^2$ is the quotient of two codimension-1 faces, say $\tau_1$ of $\sigma_i^3$ and $\delta_1$ of $\sigma^3_j$ by the identification $\phi_1:\tau_1 \to \delta_1.$ (It is possible that $i =j$). For each tetrahedron $\sigma_i^3$ and $\sigma_j^3,$ there are two normal discs, one normal triangle $p_{n_i}$ and one normal quadrilateral $q_{m_i}$ whose intersection with $\sigma^2$ is normal isotopic to $\alpha.$ In this notation, the normal coordinate $\overline{F}$ of $F$ satisfies the equations, one for each normal isotopy class of normal arc,
\begin{equation}\label{2.1} 
x_{n_i} + y_{m_i} = x_{n_j} + y_{m_j} 
\end{equation}
Define $C(K; \tri)$ to be the set of all vectors $(x_1, ..., x_{3t}, y_1, ..., y_{4t}) \in \RR^{7t}$ satisfying the equations (\ref{2.1}). 


\subsection{Tetrahedral and edge solutions}

An important theorem of Kang and Rubinstein \cite{KR1} produces a basis of $C(K; \tri).$ To state this theorem, some notation needs to be introduced. Suppose the standard basis of $\RR^{7t}$ is $\{ \epsilon_1, ..., \epsilon_{7t} \}.$ We identify $q_i$ with the vector $\epsilon_i$ and $p_j$ with $\epsilon_{3t+j}.$ For each tetrahedron $\sigma^3$ in $\tri$ containing the normal disc types
$p_{j_1}, p_{j_2}, p_{j_3}, p_{j_4}$ and $q_{i_1}, q_{i_2}, q_{i_3},$ Kang and Rubinstein associate the vector $W_{\sigma^3} = -q_{i_1} -q_{i_2} -q_{i_3} + p_{j_1}+p_{j_2}+p_{j_3} + p_{j_4}.$ For each edge $e \in  \tri,$ there are $k$ edges in the unidentified tetrahedra $\sigma_1, ..., \sigma_t$ in $X$ which are glued to form $e.$ We call $k$ the \it degree \rm of the edge $e$ in the triangulation.
 Let these $k$ edges be $e_1, ..., e_k$ in the $k$ tetrahedra $\sigma^3_{n_1}, ..., \sigma^3_{n_k}$
so that $e_i \subset \sigma^3_{n_i}.$ Note that the subscripts $n_i$ may not all be distinct. For each $\sigma^3_{n_i},$ there is a normal quadrilateral in $\sigma^3_{n_i}$ disjoint from $e_i,$ and there are two normal triangles which meet $e.$ Let $q_{l_i}$ and $p_{s_i},$ $p_{r_i}$ be their corresponding normal disc types in $\tri.$ We call $q_{l_i}$ a normal quadrilateral type \emph{facing} the edge $e.$ Kang and Rubinstein associate the vector $W_e =\sum_{i=1}^k (p_{s_i}+ p_{r_i} -q_{l_i}) \in \RR^{7t}$ to the edge $e$ of degree $k.$

It is easy to check that both vectors $W_{\sigma^3}$ and $W_{e}$ are solutions to the normal surface equations (\ref{2.1}), i.e.\thinspace $W_{\sigma^3}, W_{e} \in C(K; \tri).$ The following theorem is stated in \cite{KR2} with the hypothesis that the links of the vertices in $(K; \tri)$ are either all spheres or all tori or Klein bottles. However, the proof given in \cite{KR2} applies to the result as stated here.

\begin{theorem}[Kang-Rubinstein]\label{thm:K-R}
For any triangulation $\tri$ of a compact pseudo--manifold $K$ with or without boundary, let $\{\sigma_1^3, ..., \sigma^3_t\}$ be the set of all tetrahedra and $\{e_1, ..., e_n\}$ be the set of all edges in $\tri.$ Then $\{W_{\sigma_1^3}, ..., W_{\sigma_t^3}, W_{e_1}, ..., W_{e_n}\}$ is a basis of $C(K; \tri).$ In particular, the dimension of $C(K; \tri)$ is $t+n.$
\end{theorem}

The following lemma follows from the definition of the solutions $W_{\sigma^3}$ and $W_e.$
\begin{lemma}\label{lem:q-coords}
Under the same assumptions as in Theorem \ref{thm:K-R}:
\begin{enumerate}
\item if $q$ is a normal quadrilateral type facing edges $e_j$ and $e_k$ ($j$ could be $k$) inside $\sigma^3_i,$ then the normal coordinate corresponding to $q$ of the vector $\sum_{i=1}^t w_i W_{\sigma^3_i} +\sum_{j=1}^n z_j W_{e_j}$ is $-(w_i+z_j+z_k)$;
\item if $t$ is a normal triangle type meeting edges $e_j,$ $e_k$ and $e_l$ (the indices are not necessarily distinct) inside $\sigma^3_i,$ then the normal coordinate corresponding to $t$ of the vector $\sum_{i=1}^t w_i W_{\sigma^3_i} +\sum_{j=1}^n z_j W_{e_j}$ is $w_i+z_j+z_k+z_l.$
\end{enumerate}
\end{lemma}


\subsection{Generalised Euler characteristic function}

An embedded normal surface $S$ is properly embedded in $(K; \tri)$ and inherits a cell decomposition from its normal discs. Its Euler characteristic can therefore be determined from its normal discs by loading the contributions from edges and vertices proportionally onto the discs as follows.

If $e$ is an edge in $\tri,$ let $d(e)$ denote its degree, and if $D$ is a normal disc, let $b(D)$ be the number of normal arcs in $D \cap\partial K.$ If $\Delta$ is a normal triangle in $S$ meeting edges $e_0, e_1, e_2,$ then its contribution to the Euler characteristic of $S$ is taken to be 
\begin{equation*}
\chi^*(\Delta)=-\frac{1}{2}(1+b(\Delta))+ \frac{1}{d(e_0)}+ \frac{1}{d(e_1)}+ \frac{1}{d(e_2)}.
\end{equation*}
If $\square$ is a normal quadrilateral in $S$ meeting edges $e_0,...,e_3,$ then its contribution to Euler characteristic is taken to be 
\begin{equation*}
\chi^*(\square)=-\frac{1}{2}(2+b(\square))+ \frac{1}{d(e_0)}+ \frac{1}{d(e_1)}+ \frac{1}{d(e_2)}+ \frac{1}{d(e_3)}. 
\end{equation*}
It follows that every normal disc type $q_i$ and $t_j$ can be assigned a well--defined number $\chi^*(q_i)$ and $\chi^*(t_j)$ respectively, and one obtains a linear function $\chi^*$ on $\RR^{7t}$ defined by 
\begin{equation}\label{eq:def gen euler}
\chi^*(x_1,...,x_{3t},y_1,...,y_{4t})=\sum_{i=1}^{3t} x_i \chi^*(q_i) + \sum_{j=1}^{4t} y_j \chi^*(t_j)
\end{equation}
with the property that for any embedded normal surface, one has: $\chi(S)= \chi^*(\overline{S}),$ where $\chi(S)$ is the Euler characteristic of the surface $S,$ and $\overline{S}$ is its normal coordinate. Note that this is also true for immersed normal surfaces, and that if $S$ is a branched normal surface, then $\chi(S)\le \chi^*(\overline{S}),$ since a branch point is counted with multiplicity its index on the right hand side.

Direct computation reveals that for the canonical basis of $C(K; \tri),$ one has $\chi^*(W_{\sigma^3}) =1$ and $\chi^*(W_{e}) =2$ if $int(e)$ is contained in the interior of $K,$ and $\chi^*(W_{e}) =1$ if $int(e)$ is contained in the boundary of $K.$ The above discussion is summarised in the following:

\begin{proposition}\label{pro:gen euler}
The generalised Euler characteristic function $\chi^*:\RR^{7t}\to \RR$ defined by equation (\ref{eq:def gen euler}) is linear and has the property that for an embedded or immersed normal surface $S$ in $(K;\tri),$ $\chi^*(\overline{S}) = \chi(S),$ and for a branched immersed normal surface $S$ in $(K;\tri),$ $\chi^*(\overline{S}) = \chi(S) + \sum_v (k_v-1),$ where the sum is taken over all branch points, and $k_v$ is the branch index of $v.$
\end{proposition}


\subsection{Angle structures revisited}

Recall that an angle structure on a tetrahedron associates to each edge a positive real number, called angle, such that the sum of the angles associated to the three edges meeting in a vertex is $\pi$ for each vertex of the tetrahedron.  A simple calculation shows that this forces opposite edges to have the same angle. Since a quadrilateral disc separates a pair of opposite edges in a tetrahedron, we may consider an angle structure as a function assigning to each quadrilateral type an angle so that the sum of the angles at the three distinct quadrilateral types in a tetrahedron is $\pi.$ 

Suppose $\tri$ is an ideal triangulation of the interior $M$ of a compact 3-manifold with non-empty boundary, and let $\{q_1, ..., q_{3t}\}$ be the set of all normal quadrilateral types in $\tri.$

\begin{definition}\label{def:gen angle}
A function $\alpha\co\{q_1, ..., q_{3t}\} \to \RR$ is called a \emph{generalised angle structure on $(M; \tri)$} if it satisfies the following two properties:
\begin{enumerate}
\item If $\sigma^3 \in \tri$ is any tetrahedron and $q_i, q_j, q_k$ are the three distinct normal quadrilateral types contained in it, then
\begin{equation*}
   \alpha(q_i) + \alpha(q_j) + \alpha(q_k) =\pi.
\end{equation*}
\item If $e\in\tri$ is any edge and $\{q_{n_1}, ..., q_{n_k}\}$ is the set of all normal quadrilateral types facing it (with multiplicity, so $k$ is the degree of $e$), then
\begin{equation*}
\sum_{i=1}^k \alpha(q_{n_i})  =2\pi.
\end{equation*}
\end{enumerate}
\end{definition}

A generalised angle structure is called a \emph{semi-angle structure} if all of its values $\alpha(q_i)$ are non--negative (and hence contained in $[0,\pi]$). A generalised angle structure is called an \emph{angle structure} if all of its values $\alpha(q_i)$ are positive (and hence contained in $(0,\pi)$). By the observation at the beginning of this subsection, this new definition of an angle structure is equivalent to the one introduced by Casson and Rivin.


\section{Proofs of Theorems \ref{thm: generalised angle structures}, \ref{thm: semi--angle structures} and \ref{thm: angle structures}}
\label{sec:proofs}

If $x=(x_1,...,x_m) \in  \RR^m,$ we use $x \ge 0$ (respectively $x >0$) to mean that all components of $x_i$ are non--negative (respectively positive). The following duality result from linear programming is known as Farkas's lemma and can be found, for instance, in \cite{z}. In the following lemma, vectors in $\RR^k$ and $\RR^l$ are considered to be column vectors, and the transpose of a matrix $A$ is denoted by $A^T.$

\begin{lemma}[Farkas's Lemma]\label{lem:farkas}
  Let $A$ be a real $k\times l$ matrix, $b \in \RR^k,$ and $\cdot$
  denote the usual Euclidean inner product on $\RR^k.$
\begin{enumerate}
\item $\{ x\in \RR^l | Ax = b \} \neq \emptyset$ if and only if for
  all $y\in \RR^k$ such that $A^Ty=0,$ one has $y \cdot b=0.$
\item $\{ x\in \RR^l | Ax = b , x \ge 0 \} \neq \emptyset$ if and only if for
  all $y\in \RR^k$ such that $A^Ty\le 0,$ one has $y \cdot b\le0.$
\item $\{ x\in \RR^l | Ax = b , x > 0 \} \neq \emptyset$ if and only if for
  all $y\in \RR^k$ such that $A^Ty\neq 0$ and $A^Ty\le 0,$ one has $y \cdot
  b<0.$
\end{enumerate}
\end{lemma}


\subsection*{Proof of Theorem \ref{thm: generalised angle structures}}

Recall that $\tri$ is an ideal triangulation of the interior $M$ of a compact 3-manifold with non--empty boundary. Let $\tri^{(3)} =\{ \sigma_1^3, ..., \sigma^3_t\},$ $\tri^{(1)} =\{e_1, ..., e_n\}$ and $\{q_1, ..., q_{3t}\}$ be the sets of all tetrahedra, edges and normal quadrilateral types in $\tri.$ By Definition \ref{def:gen angle}, a generalised angle structure solves the system of linear equations listed in the definition. Assume that $Ax=b$ denotes the matrix form of this system. Then the dual $A^T$ has dual variables, one for each row in $A.$ Namely the variables are $(w, z)=(w_1, ..., w_t, z_1, ..., z_n)$ where $w_i$ is associated to the i-th tetrahedron $\sigma^3_i$ and $z_j$ is associated to the j-th edge $e_j.$ The dual system $A^T (w,z)^T =0$ is given by the following system of linear equations, one for each normal quadrilateral type $q,$
\begin{equation}\label{3.1}
w_i + z_j + z_k =0 
\end{equation}
whenever the i-th tetrahedron $\sigma^3_i$ contains the normal quadrilateral $q$ so that $q$ is facing two edges $e_j$ and
$e_k$ in $\sigma^3_i$ (with the possibility that $j=k$).

By part (1) of Farkas's lemma \ref{lem:farkas}, a generalised angle structure exists if and only if for all $(w,z) \in \RR^t \times \RR^n$ solving the equations (\ref{3.1}), the following holds:
\begin{equation}\label{3.2}
\pi \sum_{i=1}^t w_i + 2\pi \sum_{j=1}^n z_j = 0.
\end{equation}
Form the vector $W_{w,z} = \sum_{i=1}^t w_i W_{\sigma^3_i} + \sum_{j=1}^n z_j W_{e_j} \in C(M; \tri).$ Then the generalized Euler characteristic of $W_{w,z}$ is $\chi^*(W_{w,z}) = \sum_{i=1}^t w_i + 2\sum_{j=1}^n z_j$ by the definition of $\chi^*.$ Thus, by equation (\ref{3.2}), we have $\chi^*(W_{w,z})=0.$ On the other hand, by Lemma \ref{lem:q-coords}, $-(w_i+z_j+z_k)$ are the normal quadrilateral coordinates of $W_{w,z}.$ Using Theorem \ref{thm:K-R}, we conclude that the existence of the generalised angle structure is the same as for all solutions $s \in C(M; \tri)$ with zero quadrilateral coordinates, $\chi^*(s)=0.$

Now  each component $S$ of $\partial M$ is a normal surface whose normal coordinate has zero quadrilateral coordinates. On the other hand, it follows from results in \cite{to} (see also \cite{ti}) that all vectors in $C(M; \tri)$ with zero normal quadrilateral coordinates are linear combinations of the normal coordinates of the components of $\partial M.$ Thus, the existence of the generalised angle structure is equivalent to: for each component $S$ of $\partial M,$ one has $\chi^*(S) =0.$ By Proposition \ref{pro:gen euler}, $\chi^*(\overline{S}) =\chi(S).$ 

It remains to prove the claim regarding the dimension. If there exists a solution $a_0,$ then $GA(M; \tri)$ is the affine space $a_0 + \ker A,$ where $A$ is the matrix from above -- thought of as a linear
transformation from $\RR^{3t}$ to $\RR^{2t}.$ Thus, the dimension of its kernel can be calculated as $3t-2t+ \dim(\ker(A^T)).$ The argument in the previous paragraph shows that $\ker(A^T)$ is the set of all vectors in $C(M; \tri)$ with all normal quadrilateral coordinates equal to zero.
Thus by \cite{to, ti}, $\ker(A^T)$ has a basis consisting of the normal coordinates of the components of $\partial M.$ This shows that $\dim(\ker(A^T))=v,$ and hence $\dim(\ker(A))=t+v.$\qed


\subsection*{Proofs of Theorems \ref{thm: semi--angle structures} and \ref{thm: angle structures}}

As in the proof of Theorem \ref{thm: generalised angle structures}, the equivalence of statements $(1)$ and $(2)$ in Theorems 2 and 3 follows from Theorem \ref{thm:K-R} and Farkas's lemma. For example, the existence of a semi-angle structure is equivalent to: $\chi^*(W_{w,z}) \leq 0$ holds for all $(w,z) \in \RR^t \times \RR^n$ satisfying the condition that for each normal quadrilateral $q$ facing $e_j$ and $e_k$ inside the i-th tetrahedron $\sigma^3_i$:
\begin{equation}\label{(3.3)}
w_i+ z_j + z_k \leq 0.
\end{equation}
This is the same as: $\chi^*(s) \leq 0$ for all $s \in C(M; \tri)$ with non-negative quadrilateral coordinates.

Since $\chi^*$ is additive and vertex linking surfaces are tori or Klein bottles, applying Proposition \ref{pro:gen euler} shows that it is necessary and sufficient to verify the condition in (2) for all solutions $s \in C(M; \tri)$ with \emph{all} coordinates non--negative. Alluding to linearity again, it is necessary and sufficient to verify the condition in (2) for all all $s \in PC(M; \tri).$ This proves $(2) \Leftrightarrow (3).$

Since each solution in $PC(M; \tri)$ is a convex linear combination of the vertex solutions, it is necessary and sufficient to verify the condition for all vertex solutions satisfying the hypothesis; giving $(3) \Leftrightarrow (4).$\qed


\section{Angle structures for pseudo--manifolds}
\label{sec:pseudo--manifolds}

The above methods and results are now generalised to a larger class of generalised angle structures on 3-dimensional pseudo-manifolds. Section \ref{sec: Combinatorial angle structures for compact surfaces} contains a preliminary result concerning combinatorial angle structures on compact surfaces. The definition of a generalised angle structure and its associated area-curvature functions is given in Section \ref{sec: Wedges and angles}. In Section \ref{subsec: gen euler}, the generalised Euler characteristic is related to the area-curvature functions. The main results are stated in Section \ref{subsec: gen results}, their proofs are given in Section \ref{subsec: proofs gen} and an example in Section \ref{subsec: example}.


\subsection{Combinatorial angle structures for compact surfaces}
\label{sec: Combinatorial angle structures for compact surfaces}

The notion of a combinatorial angle structure, together with a combinatorial Gau\ss--Bonnet formula, can be found in Thurston's proof of Andreev's theorem (\cite{t}, Chapter 13). It has also appeared in work by Colin de Verdi\`ere \cite{CV} and Rivin \cite{Ri}. The starting point in this paper is the following generalisation of convex polygons in classical geometry.

\begin{definition}\label{def:angles for polygons}
Let $P$ be an n-sided polygon with vertices $\{v_1, ..., v_n\},$ where $n\ge 1.$ A \emph{combinatorial angle structure} on $P$ is a function $a: \{v_1, ..., v_n\} \to \RR.$ The value $a_i := a(v_i)$ is termed the (interior) angle at $v_i.$ The \emph{combinatorial area} (or simply area) of a combinatorial angle structure on $P,$ denoted by $A(a),$ is defined to be $\sum_{i=1}^n a_i - (n-2) \pi.$
\end{definition}

If $P$ is a convex n-sided polygon in the constant curvature space $\H^2,$ $\SS^2$ or $\E^2$ of curvature $\lambda,$ then $P$ inherits a natural combinatorial angle structure $a$ obtained by measuring the inner angles in the geometry. The Gau\ss-Bonnet theorem implies that the combinatorial area $A(a)$ equals $\lambda A',$ where $A'$ is the area of $P$ measured in the ambient geometry.

Not all combinatorial angle structures have geometric realisations. For example, if the three vertices of a triangle are assigned the values $a_1=\frac{1}{6}\pi,$ $a_2=a_3=\frac{2}{3}\pi,$ then $a_1+a_2+a_3>\pi,$ which shows that it cannot be realised in hyperbolic or Euclidean geometry. However, if it could be realised in spherical geometry, then there would also be a spherical triangle with interior angles $b_1=a_1=\frac{1}{6}\pi$ and $b_2=b_3=\pi - a_2 = \frac{1}{3}\pi,$ which is impossible since $b_1+b_2+b_3<\pi.$

To extend these notions to surfaces, define the \emph{corner} of a polygon $P$ at a vertex $v$ to be the collection of all open sets $U$ in the interior of $P$ so that the closure of $U$ contains the vertex $v.$ Of course, $v$ will be called the vertex of the corner.

Suppose $(S, \G)$ is a \emph{compact surface $S$ with a CW-decomposition $\G,$} i.e. $\G$ is a finite graph in $S$ with the property that each component of $S-\G$ is an open disc. Note that $\partial S \subset \G.$ The components of $S-\G$ are called the 2-cells or simply cells in $(S, \G).$ For each component $D$ of $S-\G,$ there exists an n-sided polygon $P_n$ and a continuous map $\phi: P_n \to S$ such that 
\begin{enumerate}
\item $\phi$ sends the interior of $P_n$ homemorphically onto $D,$ 
\item  $\phi$ sends vertices of $P_n$ to vertices in $\G,$ 
\item the restriction of $\phi$ to the interior of each edge in $P_n$ is injective. 
\end{enumerate}
The integer $n$ is uniquely determined. The cell $D$ is called an \emph{(open) n-sided polygon.} If $v$ is a vertex of $D,$ i.e. $v$ is a vertex in $\G$ contained in the closure of $D,$ a corner of $D$ at $v$ is defined to be the image of the corresponding corner in $P_n$ under $\phi.$ It follows that $D$ has $n$ corners, even though the closure of $D$ may contain fewer than $n$ vertices of $\G.$  For the component $D,$ a combinatorial angle structure is defined to be a function assiging a real number to each corner, and its area is defined in the same way as in Definition \ref{def:angles for polygons}.

\begin{definition}
Suppose $(S, \G)$ is a compact surface with a  CW-decomposition. A \emph{combinatorial angle structure} on $(S, \G)$ is a function $a:$  \{{\text all corners}\}$ \to \RR.$ If $c$ is a corner, we call $a(c)$ the (interior) angle of the corner. The \emph{combinatorial area $A(a)$} of $a$ is defined to be the sum of the combinatorial areas of all 2-cells in $(S, \G).$ If $v$ is a vertex in $\G,$ then the \emph{combinatorial curvature} of the angle structure at $v,$ denoted by $K_v,$ is defined as follows. Let $a_1, ..., a_k$ be the set of angles at all corners having $v$ as a vertex. If $v$ is in the interior of the surface, then $K_v := 2\pi -\sum_{i=1}^k a_i$; and if $v$ is in the boundary, then $K_v := \pi-\sum_{i=1}^k a_i.$ 
\end{definition}

The following formula relating area and curvature with the Euler characteristic is a direct consequence of the definitions.

\begin{proposition}[Combinatorial Gau\ss--Bonnet formula] \label{pro:gauss-bonnet}
If $a$ is a combinatorial angle structure on a compact surface with a CW-decomposition $(S, \G),$ then
\begin{equation*}
A(a) + \sum_{v \in \G^{(0)}} K_v = 2\pi \chi(S).
\end{equation*}
\end{proposition}

An interesting question is whether there is a combinatorial angle structure realising prescribed curvatures for vertices and areas for 2-cells in $(S, \G).$ The combinatorial Gau\ss--Bonnet formula gives a necessary condition for the existence, and the following result shows that it is also sufficient:

\begin{proposition}
Given a compact, connected surface with a CW-decomposition $(S, \G)$ and an assigment
of a number $K_v$ to each vertex $v$ and a number  $A(f)$ to each 2-cell $f,$ there exists a combinatorial angle structure on $(S, \G)$ realising these assignments if and only if
\begin{equation}\label{eq:comb gauss-bonnet}
\sum_{f \in C(2)} A(f) + \sum_{ v\in \G^{(0)}} K_v = 2\pi \chi(S),
\end{equation}
where $C(2)$ is the set of all 2-cells.
\end{proposition}

\begin{proof}
The proof is a simple application of Farkas's lemma. Label the corners by $1, ..., n,$ and denote the angle at the i-th corner by $x_i.$ A combinatorial angle structure realising the prescribed curvatures and areas is a solution $(x_1, ..., x_n)$ of the system of linear equations of the following form. At each interior  (respectively boundary) vertex $v$ adjacent to the $i_1$-th, ..., $i_k$-th corners:
\begin{equation*}
\sum_{i_s} x_{i_s} = 2\pi- K_v \ (\text{respectively } \sum_{i_s} x_{i_s} = \pi- K_v),
\end{equation*}
and at each 2-cell $f$ with corners $j_1, ..., j_m$: 
\begin{equation*}
\sum_{j_s} x_{j_s} =  A(f)+ (m-2)\pi.
\end{equation*}
The dual system has variables $y_f$ for each 2-cell $f$ and $z_v$ for each vertex $v.$ The first part of Farkas's lemma implies that the existence of the desired combinatorial angle structure is equivalent to the following. For any vector with components $y_f, z_v,$ where $f \in C(2), v \in \G^{(0)},$ satisfying $y_f+z_v=0$
whenever there is a corner $c$ inside $f$ having $v$ as a vertex, we have 
\begin{equation}\label{eq:4.2.2}
\begin{split}
\sum_{ f \in C(2)} y_f(A(f) +&(m-2)\pi)\\
&+ \sum_{v \in \G^{(0)}\cap int(S)} z_v( 2\pi -K_v)+\sum_{v \in G^{(0)}\cap \partial S} z_v(\pi -K_v)=0.
\end{split}
\end{equation}
Since $S$ is connected, the conditions $y_f+z_v=0$ imply that  there is a constant $t$ with $y_f=t$ and $z_v=-t$ for all 2-cells $f$ and vertices $v.$ It follows that $t \in \RR$ in fact parameterises the set of all vectors to be considered. Moreover,  equation (\ref{eq:4.2.2}) reduces to the Gau\ss-Bonnet formula (\ref{eq:comb gauss-bonnet}) if $t\neq 0$ and to $0=0$ otherwise. The existence of the desired angle structure is therefore equivalent to (\ref{eq:comb gauss-bonnet}).
\end{proof}

Guo \cite{Gu} recently found a characterisation of  combinatorial angle structures arising from geometric triangulations.


\subsection{Wedges and angles}
\label{sec: Wedges and angles}

Let $(K; \tri)$ be a 3--dimensional pseudo--manifold (with or without boundary) as in Section \ref{sec:preliminaries}. Let $\sigma$ be a 3--simplex in $K.$ A quadrilateral disc separates the interior of $\sigma$ into two \emph{wedges}, and there are six wedges up to normal isotopy. Normal isotopy classes of wedges will also be referred to as wedges. A \emph{generalised angle structure} on $(K; \tri)$ is a real valued function on the set of all wedges.

Given a generalised angle structure on $(K; \tri),$ any isotopy class $t$ of normal discs inherits a well-defined (combinatorial) area $A(t)$ as defined above. If $e$ is an edge in $(K; \tri)$ whose interior is contained in the interior of $K,$ then the curvature at $e,$ $\K (e),$ is defined to be $2\pi$ minus the sum of the angles of all wedges containing it; if $e$ is contained in $\partial K,$ then $\K (e)$ is defined to be $\pi$ minus the angle sum.

As in the case of surfaces, the main question of interest is the prescribing area-curvature problem. Namely, if a real number $A(t)$ is assigned to each isotopy class of normal triangles $t,$ and a real number $\K (e)$ is assigned to each edge $e,$ is there a generalised angle structure realising these assignments as its areas and curvatures? Denote the pair of functions by $(A, \kappa).$ If such a \emph{generalised angle structure with area-curvature $(A, \K)$} exists, one can further ask whether there is an assignment such that all angles are non--negative or all angles are positive. These structures will be called \emph{semi-angle} and \emph{angle structures with area-curvature $(A, \K).$} For instance, the structures discussed in Sections \ref{sec:intro}--\ref{sec:proofs} have area-curvature $(0,0).$


\subsection{Euler characteristic and area-curvature functions}
\label{subsec: gen euler}

Given a generalised angle structure on $(K; \tri)$ with area-curvature $(A, \K),$ any embedded normal surface inherits a combinatorial angle structure with respect to its induced cell decomposition. The pair of functions $(A, \K)$ encodes all relevant information except for the areas of quadrilateral discs. This motivates the following definition.

Let $N_\square$ and $N_\Delta$ be the sets of all normal isotopy classes of normal quadrilaterals and triangles respectively. For a solution $s\in C(K, \tri)$, let $x_*(s)\co N_\square \to \RR$ be the corresponding normal quadrilateral coordinate of $s,$ and $y_*(s)\co N_\Delta \to \RR$ be the corresponding normal triangle coordinate of $s.$ On $C(K; \tri),$ define a linear function $\chi^{(A, \K)}$ as follows. If $s = \sum_{i=1}^t w_i W_{\sigma^3_i} +\sum_{j=1}^n z_j W_{e_j},$ then define:
\begin{equation}\label{eq:Euler angle}
   \chi^{(A, \K)} (s) = \frac{1}{2\pi} \big(\sum_{t \in N_\Delta} y_t(s)A(t)  +  \sum_{j=1}^n 2z_j \K(e_j)  \big).
\end{equation}
Note that if $s$ is the normal coordinate of a normal surface in $(K; \tri),$ then $2z_j$ is exactly the intersection number of the surface with the edge $e_j.$ Using this observation or a direct computation using Lemma \ref{lem:q-coords} and Proposition \ref{pro:gauss-bonnet}, one obtains:
\begin{lemma} Assume that $(K; \tri)$ admits a generalised angle structure with area-curvature $(A, \K).$
\begin{enumerate}
\item If $S$ is a normal surface in $(K; \tri)$ which contains only normal triangle discs, then $\chi(S)=\chi^{(A, \K)}(\overline{S}).$
\item If $s\in C(K; \tri),$ then 
$$\chi^{(A, \K)}(s)=\chi^*(s) + \frac{1}{2\pi} \sum_{q \in N_\square}A(q)x_q(s),$$
where $A(q)$ is the combinatorial area of $q$ induced by the generalised angle structure.
\end{enumerate}
\end{lemma}


\subsection{Necessary and sufficient conditions}
\label{subsec: gen results}

We are now able to state the main results of this paper. The following proposition implies Theorem \ref{thm: generalised angle structures} since the right hand side of equation (\ref{eq:Euler angle}) vanishes under the hypothesis of the theorem.

\begin{proposition}\label{pro:gen1}
Let $(K; \tri)$ be a 3--dimensional pseudo--manifold (with or without boundary). Then $(K; \tri)$ admits a generalised angle structure with area-curvature $(A, \K)$ if and only if for each vertex linking normal surface $S$:
\begin{equation*}
\chi (S) = \chi^{(A, \K)}(\overline{S}).
\end{equation*}
If it is non--empty, then the set of all generalised angle structures with area-curvature $(A, \K)$ is an affine space of dimension $2t-n+v,$ where $v$ is the number of vertices, $n$ the number of edges and $t$ the number of tetrahedra in $\tri.$
\end{proposition}

The following two propositions are sharp in the sense that the necessary condition is in general not sufficient and vice versa when $A\neq 0,$ as illustrated by the example given in Section \ref{subsec: example}.

\begin{proposition}[Necessary for hyperbolic/Euclidean triangles]\label{pro:gen2}
Let $(K; \tri)$ be a 3--dimensional pseudo--manifold (with or without boundary) and $(A, \K)$ be prescribing area-curvature functions with the property that $A \le 0.$ 
\begin{enumerate}
\item If $(K; \tri)$ admits a semi-angle structure with area-curvature $(A, \K),$ then $\chi^* (s) \le \chi^{(A, \K)}(s)$ for all $s \in C(K; \tri)$ with all quadrilateral coordinates non--negative. 
\item If $(K; \tri)$ admits an angle structure with area-curvature $(A, \K),$ then $\chi^* (s) < \chi^{(A, \K)}(s)$ for all $s \in C(K; \tri)$ with all quadrilateral coordinates non--negative and at least one quadrilateral coordinate positive.
\end{enumerate}
\end{proposition}

\begin{proposition}[Sufficient for Euclidean/``spherical" triangles]\label{pro:gen3}
Let $(K; \tri)$ be a 3--dimensional pseudo--manifold (with or without boundary) and $(A, \K)$ be prescribing area-curvature functions with the property that $A \ge 0.$ 
\begin{enumerate}
\item If $\chi^* (s) \le \chi^{(A, \K)}(s)$ for all $s \in C(K; \tri)$ with all quadrilateral coordinates non--negative, then $(K; \tri)$ admits a semi-angle structure with area-curvature $(A, \K).$
\item If $\chi^* (s) < \chi^{(A, \K)}(s)$ for all $s \in C(K; \tri)$ with all quadrilateral coordinates non--negative and at least one quadrilateral coordinate positive, then $(K; \tri)$ admits an angle structure with area-curvature $(A, \K).$
\end{enumerate}
\end{proposition}

The following two corollaries follow directly from the above propositions, and imply Theorems \ref{thm: semi--angle structures} and \ref{thm: angle structures}.

\begin{corollary}[Semi-angle structure with Euclidean triangles]\label{cor:gen1}
Let $(K; \tri)$ be a 3--dimensional pseudo--manifold (with or without boundary) and $(0, \K)$ be prescribing area-curvature functions. Then the following are equivalent:
\begin{enumerate}
\item $(K; \tri)$ admits a semi-angle structure with area-curvature $(0, \K),$
\item $\chi^*(s)\le \chi^{(0, \K)}(s)$ for all $s \in C(K; \tri)$ with all quadrilateral coordinates non--negative,
\item $\chi (S) = \chi^{(0, \K)}(\overline{S})$ for each vertex linking normal surface $S,$ and for all \mbox{$s \in PC(K; \tri),$} one has $\chi^*(s) \le \chi^{(0, \K)}(s),$
\item $\chi (S) = \chi^{(0,\K)}(\overline{S})$ for each vertex linking normal surface $S,$ and for all vertex solutions $V_i,$ $i=1,...,m,$ one has $\chi^*(V_i)\le \chi^{(0, \K)}(V_i).$
\end{enumerate}
\end{corollary}

\begin{corollary}[Angle structure with Euclidean triangles]\label{cor:gen2}
Let $(K; \tri)$ be a 3--dimensional pseudo--manifold (with or without boundary) and $(0, \K)$ be prescribing area-curvature functions. Then the following are equivalent:
\begin{enumerate}
\item $(K; \tri)$ admits an angle structure with area-curvature $(0, \K),$
\item $\chi^*(s)< \chi^{(0, \K)}(s)$ for all $s \in C(K; \tri)$ with all quadrilateral coordinates
  non--negative and at least one quadrilateral coordinate positive, 
\item $\chi (S) = \chi^{(0, \K)}(\overline{S})$ for each vertex linking normal surface $S,$ and for all $s \in PC(K; \tri)$ with at least one quadrilateral coordinate positive, one has $\chi^*(s) <  \chi^{(0, \K)}(s),$ 
\item $\chi (S) = \chi^{(0, \K)}(\overline{S})$ for each vertex linking normal surface $S,$ and for all vertex solutions $V_i$ with at least one quadrilateral coordinate positive, one has
$\chi^*(V_i)< \chi^{(0, \K)}(V_i).$ 
\end{enumerate}
\end{corollary}


\subsection{Proofs of Propositions \ref{pro:gen1},  \ref{pro:gen2} and  \ref{pro:gen3}}
\label{subsec: proofs gen}

\begin{wrapfigure}{o}{0pt}
\psfrag{g}{{\small $\alpha$}}
\psfrag{h}{{\small $\beta$}}
\psfrag{a}{{\small $\alpha^0$}}
\psfrag{b}{{\small $\alpha^1$}}
\psfrag{c}{{\small $\alpha^2$}}
\psfrag{d}{{\small $\alpha^3$}}
\psfrag{e}{{\small $\alpha^4$}}
\psfrag{f}{{\small $\alpha^5$}}
\psfrag{0}{{\small $t^0$}}
\psfrag{1}{{\small $t^1$}}
\psfrag{2}{{\small $t^2$}}
\psfrag{3}{{\small $t^3$}}
\psfrag{q}{{\small $\square$}}
       \includegraphics[width=3.5cm]{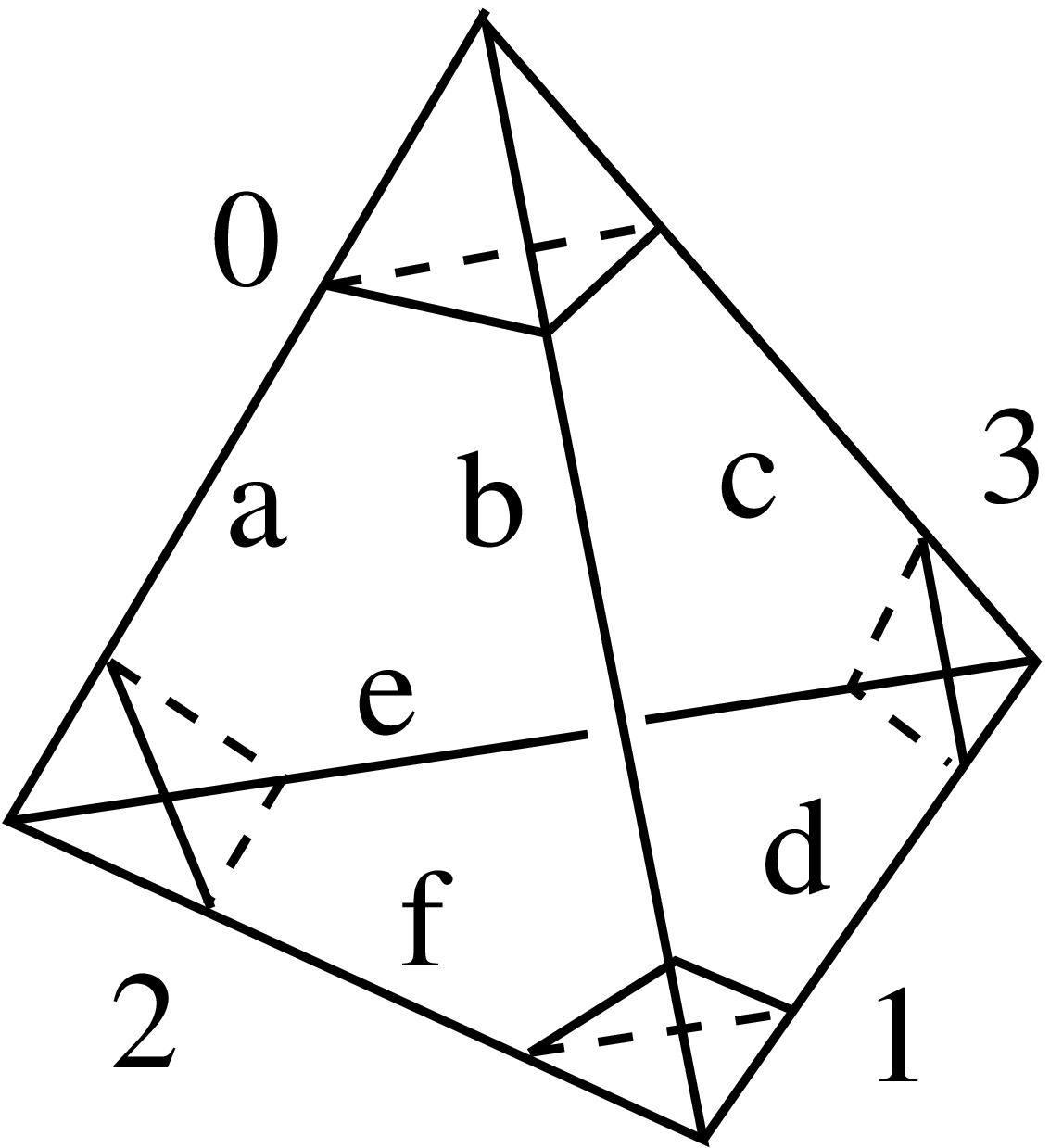}     
\end{wrapfigure}
To begin with, the problem of prescribing area and curvature is translated into the equivalent problem of prescribing angle sums for triangles and edges. The angles associated to a 3-simplex $\sigma$ are treated as variables and labelled $\alpha^0, \alpha^1, \alpha^2, \alpha^3, \alpha^4, \alpha^5,$ such that the labelling is as pictured. 

Label the 3--simplices in $K$ by $\sigma_1,...,\sigma_t,$ and their angles accordingly by $\alpha_i^k.$ For each triangle disc, the sum of its interior angles is a function of the angles associated to the wedges it is contained in:
\begin{align}
\label{eq:tri 0} \alpha_i^0 + \alpha_i^1 + \alpha_i^2 &=\t_i^0, \\
\label{eq:tri 1} \alpha_i^1 + \alpha_i^3 + \alpha_i^5 &=\t_i^1, \\
\label{eq:tri 2} \alpha_i^0 + \alpha_i^4 + \alpha_i^5 &=\t_i^2, \\
\label{eq:tri 3} \alpha_i^2 + \alpha_i^3 + \alpha_i^4 &=\t_i^3.
\end{align}

For each 1--simplex  $e_j$ in $K,$ let $a_{ij}^k \in \{0,1\}$ be the number of times $e_j$ is contained in the wedge with label $\alpha_i^k.$ Then for $i=1,...,n$ one obtains the angle sums around the edges:
\begin{equation}
\label{eq:edge j} \sum_{i=1}^t \sum_{k=0}^5 a_{ij}^k \alpha_i^k = \e_j.
\end{equation}

Given prescribing area-curvature functions $(A, \K),$ there are prescribed angle sums for triangles and edges $a_i^k = \pi + A(t_i^k)$ and $b_j = \pi -\K(e_j)$ or $b_j = 2\pi -\K(e_j)$ depending on whether $e_j$ is contained in $\partial K$ or not. Assume that $Bx=(a,b)^T$ denotes the matrix form of the system of equations given by equations (\ref{eq:tri 0}) to (\ref{eq:edge j}). Then the dual $B^T$ has one variable $h_i^k$ for each normal triangle type $t_i^k$ and one variable $z_j$ for each edge $e_j.$ The vector $B^T(h,z)^T$ has one coordinate for each wedge: if a wedge contains the edge $e_j$ and the triangles $t_i^k$ and $t_i^l,$ then the coordinate is $z_j+h_i^k + h_i^l.$

In order to link the dual $B^T$ to the normal surface theory of $(K; \tri),$ let $A$ be the matrix from the proof of Theorem \ref{thm: generalised angle structures}, thought of as arising from a system of equations for the quadrilateral types in $(K; \tri).$ Let
\begin{align*}
V=& \{ (h,z)\in\RR^{4t+n}  \mid B^T(h,z)^T\le 0 \ \text{and}\ B^T(h,z)^T\neq 0 \},\\
W=& \{ (w,z)\in\RR^{t+n} \mid A^T(w,z)^T\le 0 \ \text{and}\ A^T(w,z)^T\neq 0 \},
\end{align*}
and consider the linear transformation $\varphi\co \RR^{4t+n} \to \RR^{t+n}$ defined by $\varphi(h,z)=(w,z)$ where $w_i=h_i^0+h_i^1+h_i^2+h_i^3$ for each $i=1,...,t.$ 

\begin{lemma}\label{lem:last}
$\varphi(\ker B^T) = \ker A^T,$ $\ker B^T \cong \ker A^T$ and $\varphi(V) = W.$
\end{lemma}

\begin{proof}
If $(h,z) \in \ker B^T,$ then for each wedge, it satisfies an equation of the form:
\begin{equation*}
z_j+h_i^k + h_i^l=0.
\end{equation*}
Adding the equations for each pair of wedges sharing a quadrilateral type yields the set of equations (\ref{3.1}), which shows that $\varphi(\ker B^T) \subseteq \ker A^T.$ It also follows that the sums of coordinates of opposite pairs of edges in a tetrahedron must be equal. Using the identities $w_i=h_i^0+h_i^1+h_i^2+h_i^3$ one can then eliminate the variables associated to the triangles; if the triangle type $t_i^k$ meets edges $e_a,$ $e_b,$ $e_c,$ one obtains:
\begin{equation*}
h_i^k=-\frac{1}{2}(w_i + z_a+z_b+z_c).
\end{equation*}
This in fact yields a well-defined linear map $\ker A^T \to \ker B^T,$ which is the inverse to $\varphi.$ This shows that $\varphi|_{\ker B^T} \co \ker B^T \to \ker A^T$ is an isomorphism.

To prove the last part of the lemma, note that the definitions of $\varphi,$ $V$ and $W$ imply $\varphi(V) \subseteq W.$ Let $C_{t,n}$ be the matrix representing $\varphi,$ so $\varphi((h,z)^T) = C_{t,n} (h,z)^T.$ Let
\begin{equation*}
D_1 = 
\begin{pmatrix} 
1 & 0 & 0 & 1 & 0 & 0 \\
0 & 1 & 0 & 0 & 1 & 0 \\
0 & 0 & 1 & 0 & 0 & 1
\end{pmatrix},
\end{equation*}
and let $D_t$ be the $3t \times 6t$ block diagonal matrix with $t$ copies of $D_1$ on its diagonal. Then
\begin{equation*}
D_t B^T = A^T C_{t,n}.
\end{equation*}
Since $C_{t,n}$ has rank $t+n,$ given any $(w,z)\in \RR^{t+n}$ there exists $(h,z)\in \RR^{4t+n}$ such that $(w,z)^T= C_{t,n}(h,z)^T.$ Now $A^T(w,z)^T=A^TC_{t,n}(h,z)^T=D_tB^T(h,z)^T.$ To complete the proof of the lemma, it therefore suffices to show that for each $(h,z)^T$ such that $D_tB^T(h,z)^T\le0,$ there exists $k\in \ker C_{t,n}$ such that $B^T((h,z)^T+k)\le 0.$ This follows from the following two oberservations:

First note that with respect to a suitable labelling of the coordinates of $\RR^{6t},$ $B^T(h,z)$ lies in the intersection of the half-spaces $x_{2k-1}+x_{2k}\le 0,$ $k=1,...,3t,$ since $D_tB^T(h,z)^T\le0.$

The space $B^T(\ker C_{t,n})$ has a basis consisting of the $3t$ vectors with $(-1,1)$ in exactly one of the above coordinate pairs $(x_{2k-1}, x_{2k})$ and zeros elsewhere. This basis is obtained from the following basis for $\ker C_{t,n}.$ Note that for any element $(h,z)\in \ker C_{t,n},$ $z_1=...=z_n=0.$ The remaining coordinates come in quadruples, and a basis is given by taking all elements with all but one of the quadruples equal to zero, and the remaining equalling either $\frac{1}{2}(1,1,-1,-1),$ or
$\frac{1}{2}(1,-1,1,-1)$ or $\frac{1}{2}(1,-1,-1,1).$
\end{proof}

Given $(h,z)\in \RR^{4t+n},$ let $(w,z) = \varphi(h,z),$ and let $W_{w,z}$ denote the vector $\sum_{i=1}^t w_i W_{\sigma^3_i} + \sum_{j=1}^n z_j W_{e_j} \in C(M; \tri).$ Using Lemma \ref{lem:last} and Farkas's lemma, the proofs will be completed by looking at the following equation which follows from a direct computation:
\begin{align}\label{eq:5.0}
\notag \frac{1}{\pi} (h,z) \cdot (a,b) =& \frac{1}{\pi}  (\sum_{i=1}^n a_iz_i + \sum_{i=1}^t \sum_{k=0}^3 b_i^k h_i^k)\\
=& \chi^*(W_{w,z}) - \chi^{(A,\K)}(W_{w,z}) \\
\notag &+ \frac{1}{2\pi} \sum_{\text{wedges}} (z_j+ h_i^k + h_i^l)(b_i^k + b_i^l - 2 \pi).
\end{align}

\subsection*{Proof of Proposition \ref{pro:gen1}}
The first part of Farkas's lemma will be applied to the system $B x = (a,b)^T.$ Hence there exists a generalised angle structure with $(A, \K)$ if and only if for all $(h,z)\in \ker B^T,$ $(h,z) \cdot (a,b)=0.$ Note that 
$\frac{1}{\pi} (h,z) \cdot (a,b) = \chi^*(W_{w,z}) - \chi^{(A, \K)}(W_{w,z}) $ since for each wedge, $z_j+ h_i^k + h_i^l=0.$ Because $\varphi(\ker B^T) = \ker A^T,$ there exists a generalised angle structure with area-curvature $(A, \K)$ if and only if for all $(w,z)\in \ker A^T$:
\begin{equation*}
\chi^*(W_{w,z}) = \chi^{(A, \K)}(W_{w,z}).
\end{equation*}
As in the proof of Theorem \ref{thm: generalised angle structures}, this is the case if and only if for each $s \in C(K; \tri)$ with all quadrilateral coordinates equal to zero, $\chi^*(s) = \chi^{(A, \K)}(s).$ Since these elements are linear combinations of the normal coordinates of the vertex linking surfaces and both $\chi^*$ and $\chi^{(A, \K)}$ are linear, the equivalence stated in the theorem follows from Proposition \ref{pro:gen euler}.

The claim regarding the dimension follows from the fact that $\ker A^T \cong \ker B^T,$ so $v=\dim (\ker A^T) = \dim (\ker B^T),$ and hence $\dim (\im B) = \dim (\im B^T) = 4t+n-v,$ which implies $\dim\ker B=2t-n+v.$\qed

\subsection*{Proof of Proposition \ref{pro:gen2}}
Consider the system $B x = (a,b)^T$ when all $a_i^k\le \pi.$ The second part of Farkas's lemma implies that there is a semi-angle structure with $(A, \K)$ only if $(h,z) \cdot (a,b)\le 0$ for all $(h,z)\in V \cup \ker B^T.$ Now $B^T(h,z)^T \le 0$ implies that for each wedge, $z_j+ h_i^k + h_i^l\le 0,$ and since all $a_i^k\le \pi,$ one obtains from equation (\ref{eq:5.0}):
\begin{equation}\label{eq:6.0}
\frac{1}{\pi} (h,z) \cdot (a,b) \ge \chi^*(W_{w,z}) - \chi^{(A, \K)}(W_{w,z}).
\end{equation}
This inequality together with the fact that $\varphi(V \cup \ker B^T) = W \cup \ker A^T$ implies that there is a semi-angle structure with $(A, \K)$ only if $\chi^*(W_{w,z}) \le \chi^{(A,\K)}(W_{w,z})$ for all $(w,z)\in W \cup \ker A^T.$

The second part of the lemma follows as above from equation (\ref{eq:6.0}) and the fact that $\varphi(V) = W$ using the third part of Farkas's lemma.\qed

\subsection*{Proof of Proposition \ref{pro:gen3}}
The line of argument is similar to the previous proof. Since all $a_i^k\ge \pi,$ the inequality is reversed:
\begin{equation*}
\frac{1}{\pi} (h,z) \cdot (a,b) \le \chi^*(W_{w,z}) - \chi^{(A, \K)}(W_{w,z}),
\end{equation*}
and hence the Euler characteristic condition is observed to be sufficient for the existence of a (semi)-angle structure with $(A, \K).$\qed


\subsection{Example}
\label{subsec: example}

\begin{figure}[t]
\psfrag{a}{{$e_3$}}
\psfrag{b}{{$e_1$}}
\psfrag{c}{{$e_2$}}
\psfrag{0}{{$0$}}
\psfrag{1}{{$1$}}
\psfrag{2}{{$2$}}
\psfrag{3}{{$3$}}
\begin{center}
	\includegraphics[height=4cm]{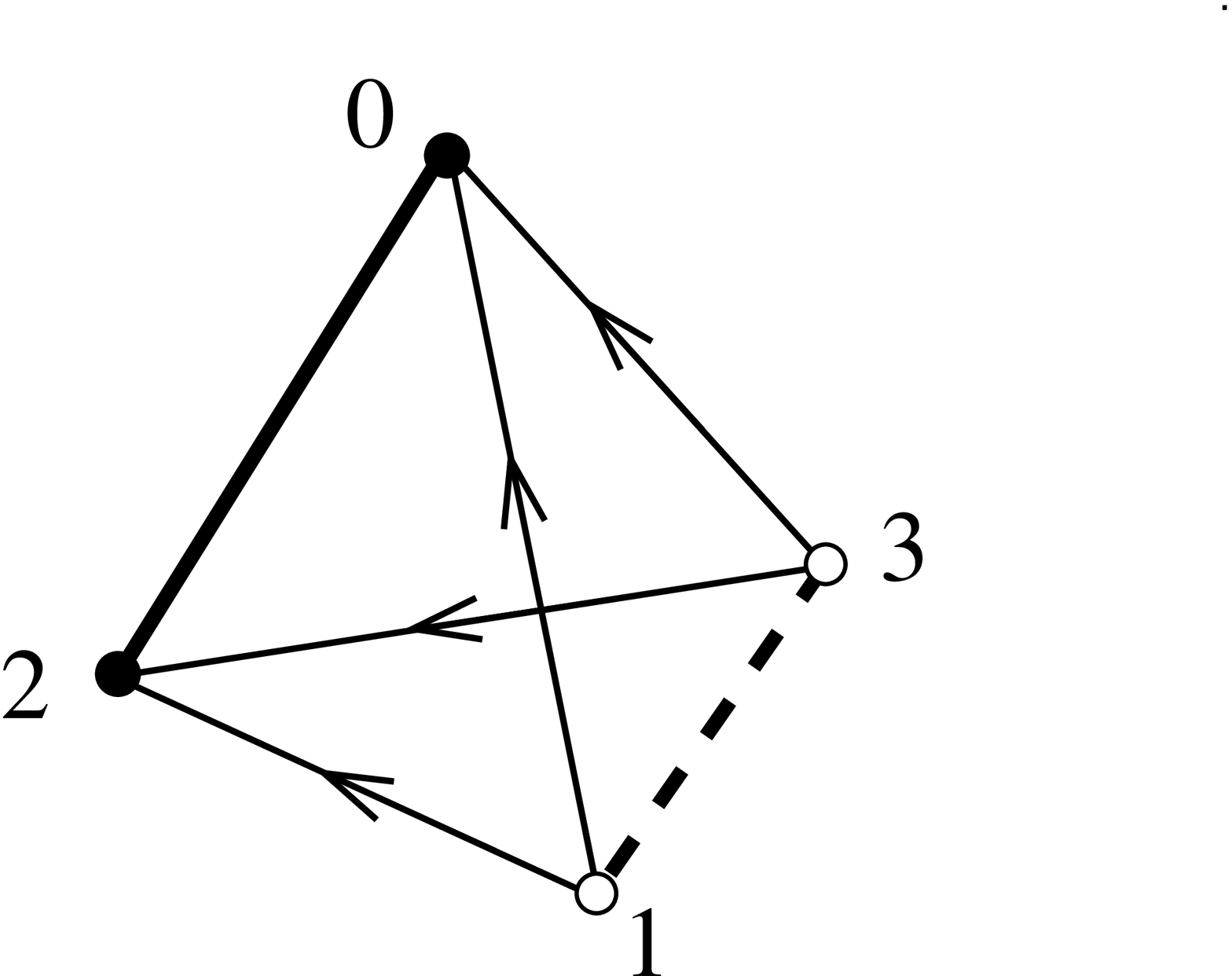}
	\quad
	\includegraphics[height=3cm]{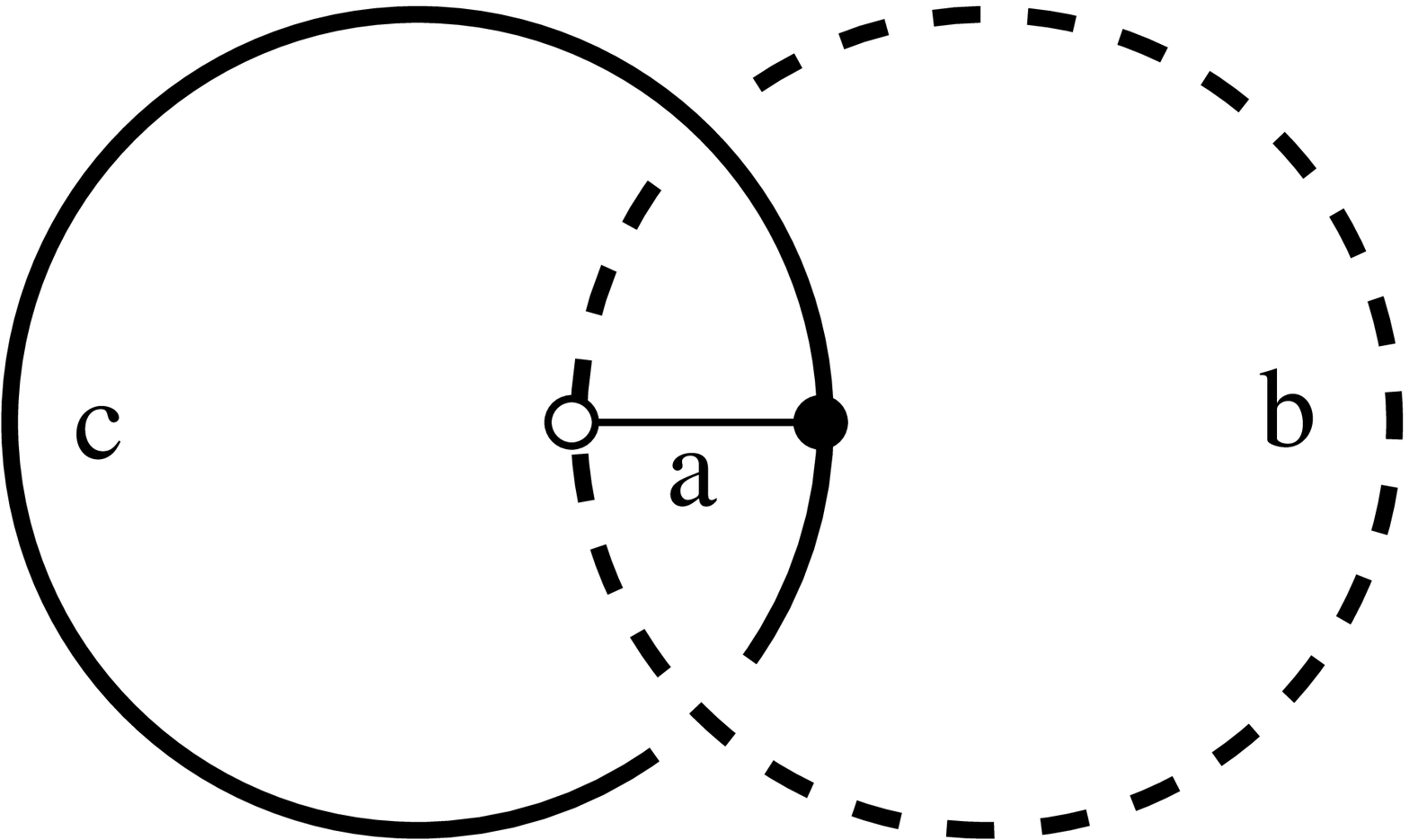}
\end{center}
    \caption{The pseudo-manifold $(K; \tri)$}
    \label{fig:example}
\end{figure}

The example given in this section is the pseudo-manifold $(K; \tri)$ obtained by taking a 3-simplex, and identifying faces in pairs such that one obtains two edges of degree one and an orientable identification space. More specifically, consider the 3-simplex $\sigma$ with labelled vertices shown in Figure \ref{fig:example}, together with the face pairings $[1,3,0]\to [1,3,2]$ and $[2,0,1]\to [2,0,3].$ Then $K$ is homeomorphic to $S^3,$ and the edges of $\tri$ form the shown graph in $S^3.$

The example is known to admit many interesting cone manifold and orbifold structures with singular locus the above graph. The following structures are described in \cite{BP}, Section 9.6. Assign the angle $\alpha$ to the two wedges containing the edges $e_1$ and $e_2,$ and the angle $\alpha/4$ to the four wedges containing the edge $e_3.$ An elementary calculation shows that there is a spherical triangle with angles $\pi, \pi/4, \pi/4$ if and only if $2\pi/3<\alpha<\pi.$ For $\alpha = 2\pi/3,$ one obtains a Euclidean triangle, for $\alpha=\pi$ a spherical bigon. Consider therefore the following cases.
\begin{enumerate}
\item[(a)] For each $\alpha \in (0, 2\pi/3),$ there is a hyperideal hyperbolic tetrahedron which realises the assignment of angles and gives $K$ minus small open neighbourhoods of the two vertices the structure of a hyperbolic cone 3-manifold with totally geodesic boundary. 
\item[(b)] For $\alpha = 2\pi/3,$ there is an ideal hyperbolic tetrahedron of shape $1+i$ realising the assignments of angles, giving $K$ minus the vertices a cusped hyperbolic 3-orbifold structure of finite volume. 
\item[(c)] If $\alpha \in (2\pi/3, \pi),$ one obtains closed cone manifolds triangulated with a single hyperbolic, Euclidean or spherical tetrahedron. The degenerate case $\alpha=\pi$ gives a spherical orbifold: it is obtained from a lens in $S^3,$ and the vertex links are spherical $(2,2,2)$--turnovers.
\end{enumerate}

Let $t^i$ denote the normal isotopy class of triangles separating the vertex $i$ of $\sigma$ from the others, and label the normal isotopy classes of quadrilateral discs as follows: $q^0$ separates the vertices in the pairs $01/23,$ $q^1$ separates $03/12,$ and $q^2$ $02/13.$ It turns out that $C(K; \tri)$ has the following basis:
\begin{equation*}
\overline{S_1}=t^0+t^2, \ \overline{S_2}=t^1+t^3, \ \overline{T}=q^2, \ \overline{R}=q^0+q^1.
\end{equation*}
The basis elements are the normal coordinates of, respectively, the vertex linking spheres $S_1$ and $S_2,$ an embedded torus $T,$ and a branched immersed real projective plane $R$ with two branch points of index two. Whence $\chi^*(\overline{S_i})=\chi(S_i)=2,$ $\chi^*(\overline{T})=\chi(\overline{T})=0$ and $\chi^*(\overline{R})=\chi(R)+2=3.$

Given the above description of a basis, one obtains that the following two conditions are equivalent:
\begin{enumerate}
\item for all $s \in C(K; \tri)$ with all quadrilateral coordinates non--negative, $\chi^*(s)\le \chi^{(A, \K)}(s),$
\item $\chi (S_i) = \chi^{(A, \K)}(\overline{S_i}),$ $\chi^*(\overline{T})\le  \chi^{(A, \K)}(\overline{T})$ and $\chi^*(\overline{R})\le  \chi^{(A, \K)}(\overline{R}),$
\end{enumerate}
and also that the following two conditions are equivalent:
\begin{enumerate}
\item for all $s \in C(K; \tri)$ with all quadrilateral coordinates non--negative and at least one quadrilateral coordinate positive, $\chi^*(s)< \chi^{(A, \K)}(s),$
\item $\chi (S_i) = \chi^{(A, \K)}(\overline{S_i}),$ $\chi^*(\overline{T}) < \chi^{(A, \K)}(\overline{T})$ and $\chi^*(\overline{R}) <  \chi^{(A, \K)}(\overline{R})$
\end{enumerate}

Letting $A_i=A(t^i)$ and $\kappa_i=\kappa(e_i),$ one obtains:
\begin{align*}
\chi^{(A, \kappa)}(\overline{S_1}) &= \frac{1}{2\pi} (A_0+A_2+2\kappa_1+\kappa_3), &
\chi^{(A, \kappa)}(\overline{T}) &=  \frac{1}{2\pi} \kappa_3, \\
\chi^{(A, \kappa)}(\overline{S_2}) &= \frac{1}{2\pi} (A_1+A_3+2\kappa_2+\kappa_3), &
\chi^{(A, \kappa)}(\overline{R}) &= \frac{1}{2\pi}(2\kappa_1+2\kappa_2+\kappa_3).
\end{align*}

It can now be verified that the structures described in (c) do not satisfy the sufficient condition given in Proposition \ref{pro:gen3} when $\alpha\ge \frac{4}{5}\pi.$ It is interesting to note that the violation is caused by the branched immersed surface, and that it occurs about $0.2$ radians after the hyperbolic to Euclidean transition. No semi-angle or angle structure can exist when $\alpha<0$; however, the necessary condition given in Proposition \ref{pro:gen2} is satisfied.


\bibliographystyle{amsplain}


\end{document}